\documentclass[12pt]{amsart}
\usepackage{amssymb}
\usepackage{amsmath}
\usepackage{amsthm}
\usepackage{times}
\usepackage{graphicx}

\usepackage[curve,all]{xy}
 \xyoption{curve}
 \xyoption{line}
\xyoption{arc}

\usepackage{color}
\usepackage{amsthm}

\newtheorem{theorem}{Theorem}[section]
\newtheorem{lemma}[theorem]{Lemma}
\newtheorem{corollary}[theorem]{Corollary}
\newtheorem{proposition}[theorem]{Proposition}

\theoremstyle{definition}
\newtheorem{definition}[theorem]{Definition}
\newtheorem{example}[theorem]{Example}

\theoremstyle{remark}
\newtheorem{remark}[theorem]{Remark}

\numberwithin{equation}{section}

\def\Ker{{\text{Ker}}}

\def\deg{{\text{deg}}}

\def\spmapright#1{\smash{%
   \mathop{\hbox to 1.3cm{\rightarrowfill}}
       \limits^{#1}}}
\def\sbmapright#1{\smash{%
   \mathop{\hbox to 1.3cm{\rightarrowfill}}
       \limits_{#1}}}
 
\newcommand{\mapright}[1]{%
\smash{\mathop{%
   \hbox to 1cm{\rightarrowfill}}\limits^{#1}}}
\newcommand{\mapleft}[1]{%
\smash{\mathop{%
   \hbox to 1cm{\leftarrowfill}}\limits^{#1}}}

\begin{document}

\title{Decomposed Richelot isogenies of Jacobian varieties  of curves of genus 3}

\author{Toshiyuki Katsura}
\address{Graduate School of Mathematical Sciences, The University of Tokyo, 
Meguro-ku, Tokyo 153-8914, Japan}
\email{tkatsura@ms.u-tokyo.ac.jp}
\thanks{Partially supported by JSPS Grant-in-Aid 
for Scientific Research (C) No. 20K03530}

\begin{abstract}
For a non-singular projective curve $C$ of genus 3 defined
over an algebraically closed field of characteristic $p \neq 2$, we give a 
necessary and sufficient condition that 
the Jacobian variety $J(C)$ has a decomposed Richelot isogeny outgoing from it
and we determine the structures of decomposed ones. 
\end{abstract}



\maketitle

\section{Introduction}
Isogeny-based cryptosystem, for example SIKE (Supersingular Isogeny Key Encapsulation),  
is one of the candidates of post-quantum cryptosystem.
The system of supersingular elliptic curves is now well examined and
achieves a great success (cf. Jao-De Feo \cite{JF} and Costello \cite{C}, for instance). 
As a next step, it is natural for researchers to investigate
higher genus cases. In the case of genus 2, 
many such trials are already done (cf. Takashima \cite{T}, Castryck--Decru--Smith \cite{CDS}
and Costello--Smith \cite{CS},
for instance) and we know now fairly well the structure of graph of
superspecial curves of genus 2 for $(2, 2)$-isogenies (cf. Ibukiyama--Katsura--Oort \cite{IKO},
Katsura--Takashima \cite{KT}, Florit--Smith \cite{FS} and Jordan--Zaytman \cite{JZ}).
As for the case of genus 3, Richelot isogenies outgoing from the Jacobian varieties of
hyperelliptic curves with tractable isotropic subgroups are studied 
(cf. Smith \cite{S}, for instance), and also Richelot isogenies outgoing from
products of 3 elliptic curves are very well analyzed (cf. Howe--Lepr\'evost- Poonen \cite{HLP}),
but general cases seem not to be well examined yet.

In this paper, we examine the decomposed Richelot isogenies outgoing
from the Jacobian varieties $J(C)$ of non-singular projective curves $C$ of genus 3 
defined
over an algebraically closed field $k$ of characteristic $p \neq 2$.
A Richelot isogeny is a $(2, 2, 2)$-isogeny outgoing
from the Jacobian variety $J(C)$ (cf. Definition \ref{Richelot isogeny}). 
Note that decomposed Richelot
isogenies (cf. Definition \ref{Richelot isogeny})
play important roles to analyze the security
of cryptosystems constructed by isogeny graph (see Costello--Smith \cite{CS}).
From a theoretical point of view, it is also interesting to examine when 
Jacobian varieties have decomposed Richelot isogenies.
In this paper, we show the following two theorems (for the definition of long automorphism
of order 2, see Section 3).

{\bf Theorem I}. 
Let $C$ be a non-singular projective curve of genus 3.
Then, there exists a decomposed Richelot isogeny outgoing from $J(C)$ 
if and only if $C$ has a long automorphism of order 2.

{\bf Theorem II}. 
Let $C$ be a non-singular projective curve of genus 3 
with a long automorphism $\sigma$ of order 2.
We set $E = C/\langle \sigma \rangle$. Then, $E$ is an elliptic curve. 
Let
$f : C \longrightarrow C/\langle \sigma \rangle = E$ be the quotient morphism,
$N_{f} : J(C) \longrightarrow E$ be the induced homomorphism and
$f^{*} :  E \cong J(E) \longrightarrow J(C)$ be the pull-back by $f$.

\begin{itemize} 
\item[$({\rm 1})$] If $C$ is hyperelliptic with hyperelliptic involution $\iota$, 
then $\{E, C/\langle \sigma\circ \iota\rangle\}$
is a set of an elliptic curve and a curve of genus 2. 
The target of the decomposed Richelot isogeny outgoing from $J(C)$ related to $\sigma$
is isomorphic to $J(E)\times J(C/\langle \sigma\circ \iota\rangle)$,
the product of Jacobian varieties. 
\item[$({\rm 2})$] If $C$ is non-hyperelliptic, then
$f^{*}$ is injective.
Moreover, $A =\Ker ~N_{f}$ is an irreducible abelian surface, and
there exist three \'etale coverings $\tilde{A}$ of $A$ of degree 2
such that the targets of the decomposed Richelot isogenies outgoing from $J(C)$ 
related to $\sigma$ are isomorphic to $(E, O) \times (\tilde{A}, \Xi)$.
Here, $\Xi$ is a principal polarization on $\tilde{A}$.
\item[$({\rm 3})$] If the Jacobian variety 
$J(C)$ has a completely decomposed Richelot isogeny, 
then $C$ is a Howe curve of genus 3. The automorphism group ${\rm Aut}(C)$ of $C$ 
contains a subgroup $G \cong {\bf Z}/2{\bf Z}\times {\bf Z}/2{\bf Z}$
with generators $\sigma$, $\tau$ such that 
the three curves $C/\langle \sigma \rangle$, $C/\langle \tau \rangle$ and
$C/\langle \sigma\circ \tau \rangle$ are elliptic curves, and 
the target of the completely decomposed Richelot isogeny
outgoing from $J(C)$ related to $\sigma$ and $\tau$ is isomorphic to
$(C/\langle \sigma \rangle, O) \times (C/\langle \tau \rangle, O) \times 
(C/\langle \sigma\circ \tau \rangle, O)$. 
\end{itemize}

We explain the outline of our paper.
In this paper, the genus of curves $C$ which we treat is always 3, 
if otherwise mentioned.
In Section 2, we prepare some lemmas which we use later. These lemmas 
are somehow known (cf. Birkenhake--Lange \cite{BL}, for instance), but
to explain our setting precisely, we give full proofs for them.
In Section 3, we examine the relation between long automorphisms of order 2
and decomposed Richelot isogenies. In Section 4, we treat the case of hyperelliptic
curves of genus 3, and give a criterion for the existence of decomposed Richelot isogenies.
In Section 5, we prepare some results 
on automorphisms of non-hyperelliptic curves of genus 3.
In Section 6, we examine the structure of Howe curves of genus 3 and show that
the Jacobian varieties of Howe curves have completely decomposed Richelot isogenies. 
This part is essentially
known in Howe--Lepr\'evost--Poonen \cite{HLP} from the dual view point of ours.
In Section 7, we treat non-hyperelliptic
curves of genus 3 and we show how non-hyperelliptic curves of genus 3 with long automorphism 
of order 2 make decomposed Richelot isogenies. Finally, summarizing our results, 
we prove Theorems I and II. Similar decompositions of Jacobian varieties are 
also investigated in Lombardo-Garc\'ia-Ritzenthaler-Sijsling \cite{LGRS}. The difference is
that their results are up to isogeny and our results are up to isomorphism.

The author thanks Katsuyuki Takashima for useful discussions and
for giving the author much information on cryptography, and
Everett Howe for useful comments 
and the information on the paper \cite{HLP}. He also thanks
the referee for his careful reading and for giving the author 
many advices.

\subsection*{Notation and conventions}

For an abelian variety $A$ and divisors $D$, $D'$ on $A$, we use the following notation.

$O$ : the zero point of $A$,

$id_A$ :  the identity of $A$,

$\iota_A$ : the inversion of $A$, i.e. the multiplication-by-(-1),

$\hat{A}= {\rm Pic}^0(A)$ : the dual (Picard variety) of $A$, 

${\rm NS}(A)$ : the N\'eron-Severi group of $A$,

$D\sim D'$: linear equivalence,

$D\approx D'$: algebraic equivalence.

For a vector space $V$ and a group $G$ which acts on $V$, we denote by $V^G$ 
the invariant subspace of $V$.
Sometimes, a Cartier divisor and the associated invertible sheaf will be identified.

\section{Preliminaries}
Let $k$ be an algebraically closed field of characteristic $p \neq 2$.
In this section, we introduce some notation and lemmas 
which we will use later.

For an abelian variety $A$ and a divisor $D$ on $A$,
we have a homomorphism
$$
\begin{array}{rccc}
\Phi_{D} :& A &\longrightarrow & {\rm Pic}^0(A) = \hat{A}\\
     & x   &\mapsto & T_x^{*}D - D
\end{array}
$$
(cf. Mumford \cite{M}). Here, $T_x$ is the translation by $x \in A$.
We put $K(D) = \Ker~\Phi_{D}$.
We know that $\Phi_{D}$ is an isogeny if $D$ is ample.
      
Let $C$ be a non-singular projective curve 
of genus $g \geq 1 $ defined over $k$. We denote by $J(C)$
the Jacobian variety of $C$, and by $\Theta$ the principal
polarization on $J(C)$ given also by $C$.
We have a natural immersion (up to translation)
$$
\alpha_C  : C \hookrightarrow J(C)= {\rm Pic}^0(C)
$$
By the abuse of terminology, we sometimes denote $\alpha_C(C)$ by $C$.
The morphism $\alpha_C$ induces a homomorphism
$$
\alpha_C^{*} : \hat{J(C)} = {\rm Pic}^0(J(C))\longrightarrow {\rm Pic}^0(C)=J(C).
$$
\begin{lemma}[Birkenhake--Lange\cite{BL}]\label{Theta}
   $\alpha_C^{*} = - \Phi_{\Theta}^{-1}$.
\end{lemma}
\proof{
We identify the image $\alpha_C(C)$ with $C$. 
As is well-known, we have $(\Theta\cdot C) = g$.
Therefore, the invertible sheaf ${\mathcal O}_{J(C)}(\Theta)\vert_{C}$
gives a divisor $\delta$ of degree $g$ on $C$.
For $x \in J(C)$, $x$ is an invertible sheaf on $C$, which we denote by ${\mathcal L}$.
Then, we have
$$
\begin{array}{rl}
\alpha_C^{*}(\Phi_{\Theta}(x)) &= (T_x^{*}(\Theta) - \Theta)\vert_C\\
   & = {\mathcal L}^{-1} \otimes \delta \otimes \delta^{-1}\\
    & = {\mathcal L}^{-1} = -x
\end{array}
$$
This means $\alpha_C^{*}\circ \Phi_{\Theta} = -{\rm id}_{J(C)}$.
Since $\Theta$ is a principal polarization, $\Phi_{\Theta}$ is
an isomorphism. Therefore, we have $\alpha_C^{*} = - \Phi_{\Theta}^{-1}$.
\hfill $\Box$}

Let $f : C \longrightarrow C'$ be a morphism of degree $2$
from $C$ to a non-singular projective curve $C'$ of genus $g'\geq 1$.
For an invertible sheaf 
${\mathcal O}_C(\sum m_iP_i) \in J(C)$ ($P_i\in C$, $m_i \in {\bf Z}$), 
the homomorphism $N_f : J(C) \longrightarrow J(C')$ is defined by 
$$
N_f({\mathcal O}_C(\sum m_iP_i)) = {\mathcal O}_{C'}(\sum m_if(P_i)).
$$
Then, by suitable choices of $\alpha_C $ and $\alpha_{C'}$, we have a commutative
diagram
$$
\begin{array}{ccc}
 C &\stackrel{\alpha_C}{\hookrightarrow} & J(C) \\
            f \downarrow \quad & & \quad \downarrow N_f \\
 C' & \stackrel{\alpha_{C'}}{\hookrightarrow} & J(C').  
\end{array}
$$

\begin{lemma}\label{N}
  $\Phi_{\Theta} \circ f^{*} = \hat{N_f}\circ\Phi_{\Theta'}$
\end{lemma}
\proof{
We have a diagram
$$
\begin{array}{ccc}
   J(C) & \stackrel{\alpha_C^{*}}{\longleftarrow} & \hat{J(C)} \\
           f^{*} \uparrow \quad & & \quad \uparrow N_f^{*} \\
 J(C') & \stackrel{\alpha_{C'}^{*}}{\longleftarrow} & \hat{J(C')}.  
\end{array}
$$ 
Therefore, using Lemma \ref{Theta}, we have
$$
\begin{array}{ccc}
   J(C) & \stackrel{\Phi_{\Theta}}{\longrightarrow} & \hat{J(C)} \\
           f^{*} \uparrow \quad & & \quad \uparrow N_f^{*} \\
 J(C') & \stackrel{\Phi_{\Theta'}}{\longrightarrow} & \hat{J(C')}.  
\end{array}
$$ 
Therefore, we have
$\Phi_{\Theta} \circ f^{*} = N_f^{*}\circ\Phi_{\Theta'}$.
Since $N_f^{*} = \hat{N}_f$, we complete our proof.
\hfill $\Box$}
\begin{lemma}\label{2Theta}
$(f^{*})^*(\Theta) \approx 2\Theta'$.
\end{lemma}
\proof{
By definition, we have $N_f\circ f^{*}=[2]_{J(C')}$. Therefore,
we have $\hat{f^{*}}\circ \hat{N}_f =[2]_{\hat{J(C')}}$.
Using Lemma \ref{N}, we have
$$
\begin{array}{rl}
   \Phi_{2\Theta'} &= [2]_{\hat{J(C')}}\circ \Phi_{\Theta'}\\
                   &=\hat{f^{*}}\circ \hat{N}_f\circ \Phi_{\Theta'}\\
                   &= \hat{f^{*}}\circ \Phi_{\Theta} \circ f^{*}\\
                   &= \Phi_{(f^*)^*(\Theta)}.
\end{array}
$$
Therefore, we have $(f^{*})^{*}(\Theta) \approx 2\Theta'$.
\hfill $\Box$}

\begin{definition}
Let $A_i$  be abelian varieties
with principal polarizations $\Theta_i$ ($i = 1, 2, \ldots, n$), respectively.
The product $(A_1, \Theta_1)\times (A_2, \Theta_2)\times \ldots \times (A_n, \Theta_n)$
means the principally polarized abelian variety $A_1 \times A_2 \times \ldots \times A_n$
with principal polarization 
$$
\Theta_1\times A_2 \times A_3 \times \ldots \times A_n
+ A_1 \times\Theta_2\times  A_3 \times \ldots \times A_n + \ldots +  A_1 \times A_2 \times \ldots \times A_{n-1} \times \Theta_n.
$$
\end{definition}

\begin{definition}\label{Richelot isogeny}
Let $C$ be a non-singular projective curve of genus $g \geq 2$,
and $J(C)$ be the Jacobian variety of $C$. We denote by $\Theta$
the canonical principal polarization of $J(C)$.
Let $A$ be an abelian variety of dimension $g$ 
with principal polarization $D$,
and $f : J(C) \longrightarrow A$ be an isogeny.
The isogeny $f$ is called a {\it Richelot isogeny} 
if $2\Theta \approx f^*(D)$. A Richelot isogeny $f$ is
said to be {\it decomposed} if there exist two principally polarized abelian varieties
$(A_i, \Theta_i)$ $(i = 1, 2)$
such that $(A, D) \cong (A_1, \Theta_1)\times (A_2, \Theta_2)$.
A decomposed Richelot isogeny is said to be 
{\it completely decomposed} if there exist elliptic curves
$E_i$ with zero point $O_i$ ($i = 1, 2, \ldots, g$)
such that 
$(A, D) \cong (E_1, O_1)\times (E_2, O_2)\times \ldots \times (E_g, O_g)$.
\end{definition}

\section{Some lemmas on automorphisms}
\begin{lemma}\label{1}
Let $C$ be a non-singular projective curve of genus $g\geq 2$,
and $\sigma$ be an automorphism of $C$ of order $n < \infty$
such that the induced
automorphism on ${\rm H}^0(C, \Omega_C^1)$ is trivial.
Then, $\sigma$ is the identity morphism.
\end{lemma}
\proof{We have a morphism $f : C \longrightarrow C/\langle \sigma \rangle$
of degree $n$. Since the induced action $\sigma^{*}$ of $\sigma$
on ${\rm H}^0(C, \Omega_C^1)$ is trivial, we have
$$
{\rm H}^0(C, \Omega_C^1) = {\rm H}^0(C, \Omega_C^1)^{\langle \sigma^{*} \rangle}
\cong {\rm H}^0(C/\langle \sigma\rangle, \Omega^1_{C/\langle \sigma \rangle})
$$
Therefore, the genus of $C/\langle \sigma\rangle$ is equal to $g$.
By the Hurwitz formula, we have
$2(g - 1) = 2n(g - 1) + \delta$ with an integer $\delta \geq 0$.
Therefore, we have $n = 1$ and $\delta = 0$. This means 
$\sigma$ is the identity morphism.
\hfill $\Box$}

\begin{lemma}\label{-1}
Let $C$ be a non-singular projective curve of genus $g\geq 3$.
If $C$ has an automorphism $\sigma$ of order 2 such that the induced
automorphism on ${\rm H}^0(C, \Omega_C^1)$ is the multiplication by $-1$,
then $C$ is a hyperelliptic curve and $\sigma$ is 
the hyperelliptic involution.
\end{lemma}
\proof{
Since ${\rm H}^0(C/\langle \sigma\rangle, \Omega^1_{C/\langle \sigma \rangle})
\cong {\rm H}^0(C, \Omega_C^1)^{\langle \sigma^{*} \rangle}
=\{0\}$, we see that the genus of $C/\langle \sigma\rangle$ is 0.
Therefore, we have the morphism 
$C \longrightarrow C/\langle \sigma\rangle\cong {\bf P}^{1}$ of degree 2.
Therefore, $C$ is hyperelliptic and $\sigma$ is the hyperelliptic involution.
\hfill $\Box$}

\begin{lemma}\label{decompose} 
Let $A$, $A_1$ and $A_2$ be abelian varieties, and let
$f : A_1 \times A_2 \longrightarrow A $ be an isogeny. 
Let $\sigma$ be an automorphism of $A$ such that 
$\sigma \circ f = f \circ(id_{A_1} \times \iota_{A_2})$
and $\Theta$ be a polarization of $A$ such that 
$\sigma^*\Theta \approx \Theta$. Then,
$$
 (A_1 \times A_2, f^*\Theta) 
 \cong (A_1, f\vert^*_{A_1}\Theta)\times (A_2, f\vert^*_{A_2}\Theta).
$$
\end{lemma}
\proof{Since $\sigma^* \Theta \approx \Theta$, we have 
$$
(id_{A_{1}}\times \iota_{A_{2}})^*(f^* \Theta) \approx (f^* \Theta).
$$
Therefore, we have 
$\Phi_{(id_{A_{1}}\times \iota_{A_{2}})^*(f^* \Theta)}
= \Phi_{f^* \Theta}$ and we have a commutative diagram
\begin{equation}\label{commutative1}
\begin{array}{rcl}
   A_1 \times A_2 & \stackrel{\Phi_{f^* \Theta}}{\longrightarrow}& 
     \hat{A}_1 \times \hat{A}_2 \\
    id_{A_1}\times \iota_{A_2} \downarrow  &    &  
    \uparrow \hat{id}_{\hat{A}_1}\times \hat{\iota}_{\hat{A}_2} \\ 
     A_1 \times A_2 & \stackrel{\Phi_{f^* \Theta}}{\longrightarrow}& 
      \hat{A}_1 \times \hat{A}_2
\end{array}
\end{equation}
We express $\Phi_{f^* \Theta}$ as a matrix
$$
\left(
\begin{array}{cc}
      \varphi_{1} & \varphi_{2} \\
      \varphi_{3} & \varphi_{4}
\end{array}
\right)
$$
(where $\varphi_{1} \in {\rm Hom}(A_1, \hat{A}_1)$, $\varphi_{2} \in {\rm Hom}(A_2, \hat{A}_1)$, 
$\varphi_{3} \in {\rm Hom}(A_1, \hat{A}_2)$ and $\varphi_{2} \in {\rm Hom}(A_2, \hat{A}_2)$).
Then, the diagram (\ref{commutative1}) says
$$
\left(
\begin{array}{cc}
      1 & 0 \\
      0 & -1
\end{array}
\right)
\left(
\begin{array}{cc}
      \varphi_{1} & \varphi_{2} \\
      \varphi_{3} & \varphi_{4}
\end{array}
\right)
\left(
\begin{array}{cc}
      1 & 0 \\
      0 & -1
\end{array}
\right)
= 
\left(
\begin{array}{cc}
      \varphi_{1} & \varphi_{2} \\
      \varphi_{3} & \varphi_{4}
\end{array}
\right).
$$
Therefore, we have $- \varphi_{2}= \varphi_{2}$ and $- \varphi_{3}= \varphi_{3}$.
Hence, we have $\varphi_{2} = 0$ and $\varphi_{3} = 0$.
This means 
$\Phi_{f^* \Theta} = \Phi_{f\vert^*_{A_1}\Theta} \times \Phi_{f\vert^*_{A_2}\Theta}$,
and we complete our proof.
\hfill $\Box$}

\begin{definition}
Let $C$ be a non-singular projective curve
of genus $g \geq 2$
and $\sigma$ be an automorphism of $C$ of order 2.
The automorphism $\sigma$ of $C$ is said to be a {\it long automorphism} 
if the $g$ eigenvalues of
the induced action of $\sigma$ on ${\rm H}^0(C, \Omega_C^1)$ are given by
$1, -1, -1, \cdots, -1$ (the number of $-1$ is $g-1$).
\end{definition}

\begin{remark}
In case $C$ is a non-singular projective curve of genus 2, 
this definition of long automorphism coincides 
with the definition of the long element
in Katsura--Takashima \cite{KT} (see also Ibukiyama-Katsura-Oort \cite{IKO}).
\end{remark}

\begin{definition}
For a polarized abelian variety with polarization $\Theta$,
we denote by ${\rm Aut}(A, \Theta)$ the group of automorphisms of $A$
which preserve the polarization $\Theta$.
\end{definition}

\begin{lemma}\label{decomposed-1} 
Let $C$ be a non-singular projective curve of genus $g \geq 2$,
and $(J(C), \Theta)$ is the Jacobian variety of $C$ with
the canonical principal polarization $\Theta$.
If the Jacobian variety $J(C)$ of $C$ has a decomposed Richelot isogeny
outgoing from $J(C)$,
then there exists an automorphism of order 2 
in ${\rm Aut}(J(C), \Theta)$ which is not the inversion. 
\end{lemma}
\proof{By assumption, we have a Richelot isogeny
\begin{equation}\label{quotient}
     \pi : J(C) \longrightarrow J(C)/G
\end{equation}
such that $G$ is a maximal isotropic subgroup of $J(C)[2]$ 
with respect to $2\Theta$, 
and that $J(C)/G$ has a decomposed principal polarization $\Theta'$
with $\pi^*\Theta' = 2\Theta$.
This means that there exist two principally polarized
abelian varieties $(A_1, \Theta_1)$ and $(A_2, \Theta_2)$
such that 
$(J(C)/G, \Theta') \cong (A_1, \Theta_1)\times (A_2, \Theta_2)$.
Since $\Theta$ is a principal polarization, we have 
an isomorphism $\varphi_{\Theta} : J(C) \cong  \hat{J}(C)$. 
By a similar reason, we have $J(C)/G \cong \hat{(J(C)/G)}$.
Using these isomorphisms, we identifies $J(C)$ (resp. $J(C)/G$)
with $\hat{J}(C)$ (resp. $\hat{(J(C)/G)}$).
Dualizing (\ref{quotient}), we have 
$$
       \eta = \hat{\pi} : J(C)/G \longrightarrow J(C).
$$
Here, we have  $J(C)/G \cong A_1 \times A_2$ with 
principal polarization $\Theta'$ such that $\eta^{*}(\Theta) \sim 2\Theta'$. 
The kernel $\Ker~ \eta$ is an isotropic subgroup of $(A_1 \times A_2)[2]$ 
with respect to the divisor $2\Theta'$. 

Since $(A_2, \Theta_2)$ is a principally
polarized abelian variety, we may assume  (by a suitable translation of $\Theta_2$)
$\iota_{A_2}^*(\Theta_2) = \Theta_2$. We set
$$
 \bar{\tau} = id_{A_1} \times \iota_{A_2} .
$$
Then, $\bar{\tau}$ is an automorphism of order 2 which is not the inversion of 
$A_1 \times A_2$. By the definition, we have
$$
\bar{\tau}^{*}(\Theta') = \Theta'.
$$
Moreover, since $\Ker~ \eta$ consists of elements of order 2 and 
$\bar{\tau}$ fixes the elements of order 2, 
$\bar{\tau}$ preserves $\Ker~ \eta$. 
Therefore, $\bar{\tau}$ induces an automorphism 
$\tau$ of $J(C) \cong (J(C)/G)/\Ker~ \eta \cong (A_1\times A_2)/\Ker~ \eta$.
Therefore, we have the following diagram:
$$
\begin{array}{ccc}
 A_1\times A_2 & \stackrel{\bar{\tau}}{\longrightarrow}  & A_1\times A_2\\
          \eta \downarrow  &        &  \downarrow \eta \\
          J(C) & \stackrel{\tau}{\longrightarrow}& J(C).
\end{array}
$$
We have
$$
  \eta^{*}\tau^{*}\Theta= \bar{\tau}^{*}\eta^{*}\Theta \sim \bar{\tau}^{*}(2\Theta') 
  = 2\Theta'.
$$
On the other hand, we have
$$
   \eta^{*}\Theta \sim 2\Theta'.
$$
Since $\eta^{*}$ is an injective homomorphism from ${\rm NS}(J(C))$ to
${\rm NS}(A_1\times A_2)$, we have $\Theta \approx \tau^{*}\Theta$. 
Therefore, $\tau$ is an element of order 2 of the group ${\rm Aut}(J(C), \Theta)$.
By definition, this is not the inversion $\iota$ of $J(C)$.
\hfill $\Box$} 
\vspace{12pt}

\section{Hyperelliptic curves of genus 3}
In this section, we assume that $C$ is a hyperelliptic curve
of genus 3. For the Jacobian variety $J(C)$ of $C$,
we denote by $\Theta$ the canonical principal polarization of $J(C)$.

\begin{proposition}\label{decomposed} 
If the Jacobian variety $J(C)$ of $C$ has a decomposed Richelot isogeny
outgoing from $J(C)$,
then there exists a long automorphism of order 2 of $C$.
\end{proposition}
\proof{
In the proof of Lemma \ref{decomposed-1}, we can 
take $A_1$ as an elliptic curve and $A_2$ as an abelian surface.
We take an automorphism $\bar{\tau} = id_{A_1} \times \iota_{A_2}$.
Then, by Lemma \ref{decomposed-1}, we have a long automorphism $\tau$
of order 2 of $J(C)$ which preserves the polarization $\Theta$.
 For hyperelliptic curves, we have 
${\rm Aut}(C)\cong {\rm Aut}(J(C), \Theta)$, and we have
${\rm H}^0(C, \Omega_{C}^{1})\cong {\rm H}^0(J(C), \Omega_{J(C)}^{1})$ 
with the compatible action of
the group of automorphisms (see Milne \cite{Mil}).
Hence, $\tau$ gives a long automorphism of order 2 of $C$.
\hfill $\Box$}

Let $\sigma$ be a long automorphism of order 2 of 
a hyperelliptic curve $C$ of genus 3, and $\iota$
be a hyperelliptic inversion of $C$. We set $\tau = \sigma\circ \iota$.
We have a morphism $\varphi : C \longrightarrow {\bf P}^1 \cong C/\langle \iota \rangle$, and
the automorphism $\sigma$ induces an automorphism of ${\bf P}^1$.
If $\sigma$ has a fixed point in the ramification points of
$\varphi$, by a suitable choice of the coordinate $x$ 
of ${\bf A}^1 \subset {\bf P}^1$,
we may assume that $\sigma$ has the fixed points at $x = 0$ and $\infty$,
and we may assume 
$$
\sigma :x \mapsto -x; \quad y \mapsto y.
$$
Then the ramification points are given by
$$
  0,  1, -1, \sqrt{a}, -\sqrt{a}, \sqrt{b}, -\sqrt{b}, \infty.
$$
Here, $a, b$ are mutually different 
and they are  equal to neither 0 nor 1.
The normal form of the curve $C$ is given by
$$
    y^2 = x(x^2 - 1) (x^2 - a)(x^2 - b).
$$
Then, the action of $\sigma$ on $C$ is 
$$
     x \mapsto -x, ~ y \mapsto \pm \sqrt{-1} y.
$$
Therefore, the order of $\sigma$ is 4, a contradiction.
Hence, $\sigma$ has no fixed points on the ramification points.
Therefore, the ramifications are given by 
$$
1, -1, \sqrt{a}, -\sqrt{a}, \sqrt{b}, -\sqrt{b}, \sqrt{c}, -\sqrt{c},
$$
and the normal form of the curve $C$ is given by
$$
    y^2 = (x^2 - 1) (x^2 - a)(x^2 - b)(x^2 - c).
$$
Elements $x^2$ and $y$ are invariant under $\sigma$.
We set $X = x^2$, $Y = y$.
Then, the defining equation of the curve $C/\langle \sigma \rangle$
is given by
$$
  Y^2 = (X - 1) (X - a)(X - b)(X - c).
$$
The curve $C/\langle \sigma \rangle$ is an elliptic curve.
We set $E_{\sigma} = C/\langle \sigma \rangle$.
We have a quotient morphism $f_1: C \longrightarrow E_{\sigma}$.
Elements $x^2$ and $xy$ are invariant under $\tau$. We set 
$X = x^2$, $Y = xy$. Then,  
the defining equation of the curve $C/\langle \tau \rangle$
is given by
$$
  Y^2 = X(X - 1) (X - a)(X - b)(X - c).
$$
The curve $C/\langle \tau \rangle$ is a curve of genus 2.
We set $C_{\tau} = C/\langle \tau \rangle$.
We have a quotient morphism $f_2 : C \longrightarrow C_{\tau}$.
Using these morphisms, we have a morphism
$$
   f = (f_1, f_2) :  C \longrightarrow E_{\sigma}\times C_{\tau}.
$$
The morphism $f$ induces a homomorphism
\begin{equation}\label{Nf}
 N_f = (N_{f_1}, N_{f_2}): J(C) \longrightarrow E_{\sigma}\times J(C_{\tau}).
\end{equation}
Note that 
$$
N_{f_1}\circ f_1^* =[2]_{E_{\sigma}}, \quad
N_{f_2}\circ f_2^* =[2]_{J(C_{\tau})}.
$$
By our construction, we have
$$
N_{f_1}\circ f_2^* =0, \quad
N_{f_2}\circ f_1^* =0.
$$
Therefore, we have
\begin{equation}\label{[2]}
N_{f}\circ f^* =[2]_{E_{\sigma}\times J(C_{\tau})}.
\end{equation}
Dualizing the situation (\ref{Nf}), we have
$$
   f^* : E_{\sigma}\times J(C_{\tau}) \longrightarrow J(C).
$$

\begin{theorem}\label{hyperelliptic}
Let $C$ is a hyperelliptic curve of genus 3 
with a long automorphism $\sigma$ of order 2.
Then, the isogeny $N_f : J(C) \longrightarrow E_{\sigma}\times J(C_{\tau})$
is a decomposed Richelot isogeny.
\end{theorem}
\proof{
Since $\sigma$ induces 
an isomorphism from $J(C)$ to $J(C)$ and we may assume
that this isomorphism is an automorphism of $J(C)$.
We have a commutative diagram
$$
\begin{array}{ccc}
  E_{\sigma}\times J(C_{\tau})& \stackrel{id_{E_{\sigma}}\times \iota_{J(C_{\tau})}}{\longrightarrow}  & E_{\sigma}\times J(C_{\tau})\\
   f^*  \downarrow    &        & \quad \downarrow f^*  \\
    \quad J(C) & \stackrel{\sigma}{\longrightarrow }   & J(C) \\
  N_f   \downarrow &     & \quad \downarrow N_f\\
   E_{\sigma}\times J(C_{\tau})& \stackrel{id_{E_{\sigma}}\times \iota_{J(C_{\tau})}}{\longrightarrow} & E_{\sigma}\times J(C_{\tau})
\end{array}  
$$
Since $\sigma^*(\Theta) = \Theta$, using Lemma \ref{decompose}, we have
$$
f^*(\Theta) \approx f_1^{*}(\Theta) \times 
J(C_{\tau})+ E_{\sigma} \times f_2^*(\Theta).
$$
Therefore, by lemma \ref{2Theta}, we see
$$
 f^*(\Theta) \approx 2(O \times J(C_{\tau}))+ 2(E_{\sigma} \times C_{\tau}).  
$$
Dualizing this situation, we have
$$
      N_f^*((O \times J(C_{\tau}))+ (E_{\sigma} \times C_{\tau}))\approx 2\Theta.
$$
This means that $N_f$ is a decomposed Richelot isogeny outgoing from $J(C)$.
\hfill $\Box$}

\section{Non-hyperelliptic curves}
In this section, we examine automorphisms of non-hyperelliptic curves.
\begin{lemma}\label{to genus 2}
Let $C$ be a non-hyperelliptic curve 
of genus 3. Then, there exist no surjective morphisms
from $C$ to curves of genus 2.
\end{lemma}
\proof{Let $C'$ be a non-singular projective curve of genus 2,
and let $f : C \longrightarrow C'$ be a nontrivial morphism.
We set $\deg~f = n$. Then we have $n \geq 2$.
If $n \geq 3$, by the Hurwitz formula, we have
$$
2(3 - 1) = n\cdot 2(2 - 1) + \delta
$$
with a non-negative integer $\delta$, which is impossible.
If $n = 2$, we have $\delta = 0$. Therefore, $f$ is an \'etale
covering. Therefore, there exists a non-trivial invertible sheaf ${\mathcal L}$ on $C'$
such that both ${\mathcal L}^{\otimes 2}$ and $f^*{\mathcal L}$ 
are trivial. Since $C'$ is of genus 2
and hyperelliptic, there exist two ramification points $P_1$, $P_2$
of the hyperelliptic covering over ${\bf P}^1$
such that $\mathcal{L} \cong {\mathcal O}_{C'}(P_2 - P_1)$, and we have
$f^*(\mathcal{L}) \cong {\mathcal O}_C$. This means $f^*(P_2) - f^*(P_1) \sim 0$,
that is, there exists a rational function $h$ on $C$ 
such that $(h) = f^*(P_2) - f^*(P_1)$.
Since $n = 2$, we see the degree of the pole divisor of $h$ is 2 and we have a morphism
$h : C \longrightarrow {\bf P}^1$ of degree 2, which contradicts the fact
that $C$ is non-hyperelliptic.
\hfill $\Box$}

\begin{corollary}\label{order-2}
Let $C$ be a non-hyperelliptic curve of genus 3, and $\sigma$
an automorphism of order 2. Then, the quotient curve $C/\langle \sigma \rangle$ 
is an elliptic curve.
\end{corollary}
\proof{
Since $C$ is non-hyperelliptic, the possibility of the genus of the curve 
$C'= C/\langle \sigma \rangle$ is either 1 or 2. However, 
2 is excluded by Lemma \ref{to genus 2}.
\hfill $\Box$}

We can also show the following corollary by the classification result
in Lombardo-Garc\'ia-Ritzenthaler-Sijsling \cite{LGRS}.

\begin{corollary}\label{1,-1,-1}
Let $C$ be a non-hyperelliptic curve of genus 3 and
$\sigma$ is an automorphism of $C$ of order 2 .
Then, the eigenvalues of the action of $\sigma^*$ on
${\rm H}^0(C, \Omega_C^{1})$ are $1, -1, -1$, that is,
$\sigma$ is a long automorphism.
\end{corollary}
\proof{
By Lemmas \ref{1} and \ref{-1}, we can exclude
$\{1, 1, 1\}$ and $\{-1, -1, -1\}$. 
Suppose the eigenvalues are $1, 1, -1$.
Then, we have
$$
\dim {\rm H}^0(C/\langle \sigma \rangle, \Omega^1_{C/\langle \sigma \rangle})
= \dim {\rm H}^0(C, \Omega^1_{C})^{\langle \sigma \rangle} = 2,
$$
that is, the genus of the curve $C/\langle \sigma \rangle$ is equal to 2,
which is excluded by Lemma \ref{to genus 2}.
\hfill $\Box$}

\begin{proposition}\label{decomposed2} 
Let $C$ be a non-hyperelliptic curve of genus 3.
If $C$ has a decomposed Richelot isogeny
outgoing from $J(C)$,
then there exists a long automorphism of order 2 of $C$ .
\end{proposition}
\proof{By Lemma \ref{decomposed-1}, we have a long automorphism $\tau$
of order 2
of $J(C)$ which preserves the polarization $\Theta$.
For non-hyperelliptic curves, either $\tau$ or $-\tau$
is induced from an element of
${\rm Aut}(C)$ (cf. Milne \cite{Mil}). We have an isomorphism
${\rm H}^0(C, \Omega_{C}^{1})\cong {\rm H}^0(J(C), \Omega_{J(C)}^{1})$ 
with the compatible actions of automorphisms in ${\rm Aut}(C)$.
By Corollary \ref{1,-1,-1}, $-\tau$ cannot become an automorphism of $C$.
Therefore, $\tau$ comes from an automorphism of $C$.
Hence, this gives a long automorphism of order 2 of $C$.
\hfill $\Box$}

\section{Howe curves}
Let $E_1$, $E_2$ be two elliptic curves, and
let $f_1 : E_1\longrightarrow {\bf P}^1$, $f_2 : E_2 \longrightarrow {\bf P}^1$ 
be morphisms of degree 2.
We consider the fiber product $E_1 \times_{{\bf P}^1} E_2$:
$$
\begin{array}{ccc}
   E_1 \times_{{\bf P}^1} E_2 & \stackrel{\pi_2}{\longrightarrow} & E_2 \\
     \pi_1 \downarrow &        & \downarrow  f_2\\
      E_1     & \stackrel{f_1}{\longrightarrow}   & {\bf P}^1.
\end{array}
$$
We denote by $r$ the number of common ramification points
of $f_1$ and $f_2$ ($0 \leq r < 4$). We exclude the case $r = 4$. Because
if $r = 4$, there exists an isomorphism $\varphi : E_1 \cong E_2$ such that
$f_2 \circ \varphi = f_1$ and the fiber product $E_1 \times_{{\bf P}^1} E_2$
is not irreducible.
We denote by $C$ the non-singular projective model of $E_1 \times_{{\bf P}^1} E_2$,
and  we denote by $h : C \longrightarrow E_1 \times_{{\bf P}^1} E_2$ 
the resolution of singularities.
We call $C$ a Howe curve (cf. Howe \cite{Ho} and Kudo-Harashita-Senda \cite{KHS}). 
Note that in case $C$ is a curve of genus 3, then this curve $C$ is historically
called a Ciani curve  (cf. Ciani \cite{Ci}).
There exist two automorphisms $\sigma$, $\tau$ of order 2
of $C$
such that $C/ \langle \sigma \rangle \cong E_1$ and
$C/ \langle \tau \rangle \cong E_2$. It is clear
that $\langle \sigma, \tau\rangle \cong {\bf Z}/2{\bf Z} \times {\bf Z}/2{\bf Z}$.
We set $h_1 = \pi_1\circ h$. Then, the degree of $h_1$ is 2.
The genus of a Howe curve is given by the following proposition.
\begin{proposition}\label{Howe}
The genus of $C$ is equal to $5 - r$.
\end{proposition}
\proof{
Let $P \in {\bf P}^1$ be a common ramification point of $f_1$ and $f_2$.
We can choose a coordinate $x$ on ${\bf A}^1 \subset {\bf P}^1$
such that $P$ is given by $x = 0$. Then, the equation
of $E_1$ (resp. $E_2$) around $P$ is given by
$$
   y_1^2 = u_1 x\quad (\mbox{resp.}~y_2^2 = u_2 x).
$$
Here, $u_1$ and $u_2$ are units at $P$.
We denote by $\tilde{P}$ the point 
of the fiber product $E_1 \times_{{\bf P}^1} E_2$
over $P$. Then, around $\tilde{P}$ 
the fiber product $E_1 \times_{{\bf P}^1} E_2$
is defined by
$$
   y_1^2 = u_1 x, ~y_2^2 = u_2 x.
$$
Therefore, by eliminating $x$, the equation around $\tilde{P}$
is given by the equation $u_2y_1^2 = u_1y_2^2$. This means
that $\tilde{P}$ is a singular point with two branches.
Therefore, on $C$ $\tilde{P}$ splits into two non-singular
points and  $P$ is not a ramification point of $h_1$.

By the meaning of fiber product, the branch points of $f_1$
whose images by $f_1$ are not ramification points of $f_2$
are not ramification points of $h_1$, and the points on
$E_1$ which are not 
branch points of $f_1$ and whose images by $f_1$ are 
ramification points of $f_2$ are ramification points 
of $h_1$. Therefore, 
on the curve $C$, $h_1$ has $2(4 - r)$ branch points
of index 2. Applying the Hurwitz formula 
to the morphism $h_1: C \longrightarrow E_1$, we have
$$
   2(g(C) -1) = 2\cdot 2(g(E_1) -1) + 2(4 - r)
$$
Since $g(E_1) = 1$, we have the result.
\hfill $\Box$}

The following two theorems are essentially known in
Howe--Lepr\'evost--Poonen \cite{HLP}.
\begin{theorem}\label{completely decomposed}
Let $C$ be a Howe curve of genus 3. 
Then, there exists a completely decomposed 
Richelot isogeny outgoing from $J(C)$.
\end{theorem}
\proof{
We set $E_3 = C/\langle \sigma\circ \tau\rangle$.
Since $C/ \langle \sigma \rangle \cong E_1$ (resp.
$C/ \langle \tau \rangle \cong E_2$) is
an elliptic curve,
the eigenvalues of the action of $\sigma$ (resp. $\tau$)
on ${\rm H}^0(C, \Omega_C^1)$ are given 
by $1, -1, -1$ (resp. $-1, 1, -1$)
with respect to a suitable choice of the basis of
${\rm H}^0(C, \Omega_C^1)$.
Therefore, the eigenvalues of the action of $\sigma \circ \tau$
on ${\rm H}^0(C, \Omega_C^1)$ are given by $-1, -1, 1$.
Therefore, $E_3$ is an elliptic curve.
We denote by $\Theta$ the canonical principal divisor
of $J(C)$.
By the universality of Jacobian variety 
we have an isogeny
$$
    f: J(C) \longrightarrow E_1 \times E_2 \times E_3.
$$
Then by a similar method to the one in Theorem \ref{hyperelliptic},
we have $2\Theta = f^*(\{0\} \times E_2\times E_3 + 
E_1 \times \{0\}\times E_3 + E_1 \times E_2\times \{0\})$
and $f$ is a completely decomposed Richelot isogeny.
\hfill $\Box$}

\begin{theorem}\label{Howe curve of genus 3}
Let $C$ be a non-singular curve of genus 3.
If there exists a completely decomposed Richelot isogeny
outgoing from $J(C)$, then $C$ is a Howe curve of genus 3.
\end{theorem}
\proof{As in Propositions \ref{decomposed} 
and \ref{decomposed2},
we have two automorphisms $\sigma$ , $\tau$ of $C$ of order 2
such that $\sigma\circ \tau = \tau \circ\sigma$.
The eigenvalues of the actions of $\sigma$ and $\tau$
are both given by one 1 and two $-1$'s.  
Therefore, $E_1 =C/\langle \sigma\rangle$
and $E_2 =C/\langle \tau \rangle$ are elliptic curves.
The automorphism $\tau$ (resp. $\sigma$)
induces the inversion of $E_1$ (resp. $E_2$), and
$C/\langle \sigma, \tau\rangle \cong {\bf P}^1$.
Then, considering the fiber product, 
we have a commutative diagram:
$$
\begin{array}{ccc}
   E_1 \times_{{\bf P}^1} E_2 & \longrightarrow & E_2 \\
      \downarrow &        & \downarrow  f_2\\
      E_1     & \stackrel{f_1}{\longrightarrow}   & {\bf P}^1.
\end{array}
$$
Since we have  morphisms $C \longrightarrow E_1$ and  
$C \longrightarrow E_2$, 
by the universality of fiber product there exists a morphism 
$f : C \longrightarrow  E_1 \times_{{\bf P}^1} E_2$.
By the degree calculation of morphisms, we see $\deg ~f = 1$.
Therefore, $C$ is birationally equivalent to $E_1 \times_{{\bf P}^1} E_2$
and $C$ is a Howe curve.  
\hfill $\Box$}

Many examples of Howe curves are known (cf. Howe--Lepr\'evost--Poonen \cite{HLP}
and Brock \cite{B}). We give here typical examples of
a hyperelliptic Howe curve and a non-hyperelliptic one. The examples
are well-known (cf. Lombardo-Garc\'ia-Ritzenthaler-Sijsling \cite{LGRS}),
but such simple examples make our situation clearer.

\begin{example} We consider the non-singular complete model $C$ of a curve
defined by
$$
   y^2 = x^8 -1.
$$
The genus of $C$ is 3 and it has two automorphisms defined by
$$
\sigma : x \mapsto -x,~y \mapsto y; \quad \tau: x \mapsto \zeta/x, ~y \mapsto \zeta^2y/x^4.
$$
Here, $\zeta$ is a primitive eighth root of unity. 
Then, they are long automorphisms of order 2 with $\sigma \circ \tau = \tau \circ \sigma$. 
Therefore, by the proof of Theorem \ref{Howe curve of genus 3},
$C$ is a  hyperelliptic Howe curve.
\end{example}

\begin{example} We consider the non-singular complete model $C$ of a Fermat curve
defined by
$$
   x^4 + y^4 = 1.
$$
The genus of $C$ is 3 and it has two automorphisms defined by
$$
\sigma : x \mapsto -x,~y \mapsto y; \quad \tau: x \mapsto x, ~y \mapsto -y.
$$ 
Then, they are long automorphisms of order 2 with $\sigma \circ \tau = \tau \circ \sigma$. 
Therefore, $C$ is a non-hyperelliptic Howe curve.
\end{example}

\section{Non-hyperelliptic curves of genus 3 with long automorphism}

Let $C$ be a non-hyperelliptic curve of genus 3 
with an automorphism $\sigma$ of order 2. By Corollary \ref{order-2}
the quotient curve $E = C/\langle \sigma \rangle$ is an elliptic curve,
and we have the quotient morphism $f : C \longrightarrow E$. 
As before, choosing an immersion $\alpha= \alpha_C : C \hookrightarrow J(C)$
suitably, we have a commutative diagram
$$
\begin{array}{rcl}
      C & \stackrel{\alpha}{\hookrightarrow}&  J(C)\\ 
            &f \searrow \quad &   \downarrow N_f \\
                &   &   E.
\end{array}
$$
\begin{lemma}\label{injective1}
$f^{*}: J(E) \longrightarrow J(C)$ is injective.
\end{lemma}
\proof{
Suppose that $f^{*}$ is not injective. We denote the zero element of $E$ by $O$.
Since any element of $J(E)$ is given by $P - O$ with a suitable point $P \in E$,
there exists a point $Q$ ($Q \neq O$) of $E$ such that $f^{*}(Q - O)$ 
is linearly equivalent to 0.
This means there exists a rational function $h$ on $C$ such that
$(h) = f^{*}(Q) - f^{*}(O)$. Since $f$ is degree 2, we have a morphism
$h :C \longrightarrow {\bf P}^{1}$ which is of degree 2. This contradicts
the assumption that $C$ is not hyperelliptic.
\hfill $\Box$}

We set $\Ker ~N_{f} = A$. We denote by $i_{A}$ the natural immersion of $A$
into $J(C)$:
$$
     i_{A}: A \hookrightarrow J(C).
$$
\begin{lemma}
$\alpha(C)\cdot A = 2$.
\end{lemma}
\proof{
For the zero point $O \in E$,
we have 
$$
\alpha( C)\circ A = \deg~ (\alpha^{-1}\circ N_{f}^{-1}(O)) = \deg f^{-1}(O) = 2.
$$
\hfill $\Box$}
\begin{lemma}\label{irreducible} 
$A$ is irreducible.
\end{lemma}
\proof{ 
Since $\alpha(C)\cdot A = 2$,
the curve $\alpha (C)$ will intersect $A$ with two points.
If $A$ is not irreducible, then considering the Stein factorization,
we have a fiber space such that $\alpha (C)$  is  a section of the
fiber space. However, since $J(C)$ is an abelian variety,
the base curve is an elliptic curve. Therefore, the curve of genus 3
cannot become a section.
\hfill $\Box$}

For the canonical principal polarization $\Theta$ of $J(C)$,
we set $D = A\cap \Theta$. Then, $D$ is a divisor on the abelian surface $A$.
\begin{lemma}\label{A}
$i_{A}^{*}(\Theta) = D$ and $D^{2} = 4$.
\end{lemma}
\proof{
The former part comes from the definition.
By Matsusaka's theorem on the characterization of Jacobian variety, 
we have $(1/2!)\Theta^{2}\approx \alpha(C)$. 
Therefore, we have
$$
   D^{2} =(\Theta\cdot (\Theta\cdot A)) = (\Theta^{2}\cdot A) =2(\alpha (C)\cdot A)= 4.
$$
\hfill $\Box$}

By the identification of $E$ with $\hat{E}$, we can regard
$f^{*}$ as the natural immersion $i_{E} : E \hookrightarrow J(C)$.

\begin{lemma}\label{E}
   $f^{*}(\Theta) \approx 2O$.
\end{lemma}
\proof{
This follows from Lemma \ref{2Theta}.
\hfill $\Box$}
\begin{lemma}\label{2,2}
Let $L$ be an ample divisor on an abelian surface $A$
with $\vert K(L)\vert = 4$.
Then, $K(L) \cong {\bf Z}/2{\bf Z} \times {\bf Z}/2{\bf Z}$.
\end{lemma}
\proof{
Suppose that  $K(L) \cong {\bf Z}/4{\bf Z}$. Since $e^{L}$ is alternating, 
for a generator $\zeta \in K(L)$ we have $e^{L}(\zeta, \zeta) = 1$, which
contradicts the fact that $e^{L}$ is a non-degenerate pairing on $K(L)$
(cf. Mumford \cite{M}).
\hfill $\Box$}

\begin{lemma}\label{2,2,4}
Let $L$ be an ample divisor on an abelian surface $A$.
Then, $K(L)$ cannot be isomorphic to 
${\bf Z}/2{\bf Z} \times {\bf Z}/2{\bf Z}\times {\bf Z}/4{\bf Z}$.
\end{lemma}
\proof{
Suppose that $K(L) \cong {\bf Z}/2{\bf Z} \times {\bf Z}/2{\bf Z}\times {\bf Z}/4{\bf Z}$.
Then, the subgroup $G \cong {\bf Z}/4{\bf Z}$ of $K(L)$
is an isotropic subgroup with respect to the pairing $e^{L}$ 
as in the proof of Lemma \ref{2,2}.
Therefore, we have a principal divisor $\Xi$ on $A/G$ and a commutative diagram
$$
\begin{array}{rcc}
        A   &\stackrel{\Phi_{L}}{\longrightarrow} & \hat{A}\\
       \pi \downarrow &      & \uparrow \hat{\pi} \\
       A/G & \stackrel{\Phi_{\Xi}}{\longrightarrow} & \hat{A/G}
\end{array}
$$
Note that $\Phi_{\Xi}$ is an isomorphism.
Since $K(L) \cong \Ker~ \Phi_{L}\cong {\bf Z}/2{\bf Z} \times {\bf Z}/2{\bf Z}\times {\bf Z}/4{\bf Z}$ and $\Ker~\pi \cong G \cong {\bf Z}/4{\bf Z}$, we see that 
$\Ker~\hat{\pi} \cong  {\bf Z}/2{\bf Z} \times {\bf Z}/2{\bf Z}$,
which contradicts the fact that $\Ker~\hat{\pi}$ is dual to $\Ker ~\pi$ 
(cf. Mumford \cite{M}).
\hfill $\Box$}
 
By abuse of notation, we denote by $E$ the image of $f^{*}$.
Then, we have a homomorphism
$$
     i_{E} + i_{A} : E \times A \longrightarrow J(C).
$$
\begin{lemma}\label{Phi}
$\Phi_{(i_{E} + i_{A})^* \Theta} = \Phi_{i_{E}^{*}\Theta} \times \Phi_{i_{A}^{*}\Theta}$.
\end{lemma}
\proof{
On $E$ $\sigma$ acts as the identity and on $A$ $\sigma$ acts as the inversion 
$\iota_{A}$ of $A$. Therefore, we have a commutative diagram
$$
\begin{array}{rcl}
    E \times A & \stackrel{id_{E}\times \iota_{A} }{\longrightarrow} & E \times A  \\
          i_{E} + i_{A}   \downarrow  &       & \downarrow i_{E} + i_{A}\\
          J(C) & \stackrel{\sigma}{\longrightarrow}    &   J(C) 
\end{array}
$$
and since $\sigma^*\Theta \approx \Theta$, we get our result 
by Lemma \ref{decompose}.
\hfill $\Box$} 

\begin{corollary}
$\Phi_{(i_{E} + i_{A})^*\Theta} = \Phi_{2O} \times \Phi_{D}$.
\end{corollary}
\proof{
This follows from Lemmas \ref{A}, \ref{E} and \ref{Phi}.
\hfill $\Box$}

Since $D^{2} = 4$, we have $\vert K(D) \vert = ((D)^{2}/2)^{2} = 4$.
Therefore, by Lemma \ref{2,2} we see 
$K(D) \cong {\bf Z}/2{\bf Z} \oplus {\bf Z}/2{\bf Z}$.
Therefore, we have a homomorphism $\varphi : \hat{A} \longrightarrow A$
such that $\Phi_{D}\circ \varphi = [2]_{\hat{A}}$.
Since $\Ker ~\varphi \cong {\bf Z}/2{\bf Z} \oplus {\bf Z}/2{\bf Z}$,
we have three elements of order 2 in $\Ker~ \varphi$.
We take one of them, say $a \in \Ker ~\varphi$, $a \neq 0$.
Then, we have the following homomorphisms:
$$
[2]_{\hat{A}} : \hat{A} \longrightarrow \hat{A}/\langle a\rangle \stackrel{\pi}{\longrightarrow}
A \stackrel{\Phi_{D}}{\longrightarrow} \hat{A}.
$$
We set $\tilde{A}_{a} = \hat{A}/\langle a\rangle$
Using this decomposition of the homomorphism $[2]_{\hat{A}}$, we have a diagram
$$
\begin{array}{lcl}
       \tilde{A}_{a} &\stackrel{\Phi_{\pi^{*}D}}\longrightarrow &{\hat{\tilde{A}}_{a}}\\
      \downarrow \pi&       &\uparrow \hat{\pi}\\
          A   & \stackrel{\Phi_{D}}{\longrightarrow} &  \hat{A}.
\end{array}
$$
Since $(\pi^{*}D)^{2}= (\deg~\pi) (D^{2}) = 8$, we have 
$\deg~ \Phi_{\pi^{*}D} = ((\pi^{*}D)^{2}/2)^{2} = 16$. Therefore,
we have $\vert K(\pi^{*}D)\vert = 16$. Since 
$K(\pi^{*}D) \supset \Ker~\Phi_{D}\circ\pi 
\cong {\bf Z}/2{\bf Z} \oplus {\bf Z}/2{\bf Z} \oplus {\bf Z}/2{\bf Z}$,
$K(\pi^{*}D)$ is isomorphic to either 
${\bf Z}/2{\bf Z} \oplus {\bf Z}/2{\bf Z} \oplus {\bf Z}/4{\bf Z}$
or ${\bf Z}/2{\bf Z} \oplus {\bf Z}/2{\bf Z} \oplus {\bf Z}/2{\bf Z}\oplus {\bf Z}/2{\bf Z}$.
By Lemma \ref{2,2,4}, we conclude  
$$
K(\pi^{*}D)\cong {\bf Z}/2{\bf Z} \oplus {\bf Z}/2{\bf Z} 
\oplus {\bf Z}/2{\bf Z}\oplus {\bf Z}/2{\bf Z}.
$$
Namely, we have $K(\pi^{*}D)\cong \Ker~[2]_{\tilde{A}_{a}}$.
By Mumford \cite[Section 23, Theorem 3]{M}, we see that there exists
a principal divisor $\Xi$ on $A$ such that $\pi^{*}D \approx 2\Xi$.
Hence, we have the following theorem.

\begin{theorem}\label{non-hyperelliptic}
Let $C$ be a non-hyperelliptic curve of genus 3 with an automorphism $\sigma$
of order two. Then, related to the automorphism $\sigma$,
there exist three decomposed Richelot isogenies
outgoing from the Jacobian variety $J(C)$.
\end{theorem}
\proof{
Using the notation above, we consider the isogeny
$$
     \tilde{\rho} : E \times \tilde{A}_{a} \stackrel{id_{E}\times \pi}{\longrightarrow}
     E \times A \stackrel{i_{E} + i_{A}}{\longrightarrow} J(C).
$$
Then, we have $\tilde{\rho}^{*}\Theta = 2(O\times \tilde{A }_{a}+ E \times \Xi)$.
Therefore, there exists a homomorphism $\rho :J(C) \longrightarrow E \times \tilde{A}_{a}$
such that $\rho\circ \tilde{\rho} = [2]_{E \times \tilde{A}_{a}}$ and
$\rho^{*}(O\times \tilde{A}_{a} + E \times \Xi) = 2\Theta$.  We have 3 possibilities
for the choice of $a$.
\hfill $\Box$}

\begin{remark}
The decomposition, up to isogeny, of the Jacobian variety of a curve of genus 3 
with automorphism of order 2 into a product of an elliptic curve and
the Jacobian variety of a curve of genus 2 is studied and classified 
in Lombardo, Garc\'ia, Ritzenthaler and Sijsling \cite{LGRS}.
They give the concrete equations for the curves of genus 2.
Although the subject of the paper is not Richelot isogeny,
their result is closely related to Theorem \ref{non-hyperelliptic}. 
\end{remark}

Now, we are ready to show Theorems I and II.
As for Theorem I, the result that if $C$ has a long automorphism of order 2, 
then we have a decomposed
Richelot isogeny outgoing from $J(C)$ is shown in  Theorem \ref{hyperelliptic} in case
$C$ is hyperelliptic, and in Theorem \ref{non-hyperelliptic} in case $C$ is non-hyperelliptic.
The converse follows from Proposition \ref{decomposed} in case $C$ is hyperelliptic,
and from Proposition \ref{decomposed2} in case $C$ is non-hyperelliptic.
As for Theorem II, the result in (1) on the target of the decomposed Richelot isogeny
outgoing from $J(C)$ is given in Theorem \ref{hyperelliptic}.
The injectivity of $f^{*}$ in (2) is proved in Lemma \ref{injective1}, and
the irreducibility of $A = \Ker~N_f$ is proved in Lemma \ref{irreducible}.
The structure of decomposed Richelot isogeny of the Jacobian
variety $J(C)$ of non-hyperelliptic curve $C$ of genus 3 with
an automorphism of order 2 in (2) is given in Theorem \ref{non-hyperelliptic} with Corollary \ref{order-2}, 
which is the most important result in this paper. The former part of
(3) is proved in Theorem \ref{Howe curve of genus 3} and the latter part of
(3) is proved in Theorem \ref{completely decomposed}.

Finally, we examine the number of superspecial curves of genus 3
whose Jacobian varieties $J(C)$ have decomposed Richelot isogenies outgoing from $J(C)$.
\begin{proposition} 
Let $k$ be an algebraically closed field of characteristic $p > 2$.
Then, the asymptotic behavier of  
the rate of superspecial curves $C$ of genus 3 
whose Jacobian varieties $J(C)$ have decomposed Richelot isogenies 
outgoing from $J(C)$ 
to the superspecial curves of genus 3 is given by
$$
   \frac{1260}{p^2}.
$$
\end{proposition}
\proof{
By Brock \cite{B} and Hashimoto \cite{H},
the main term of the number of superspecial curves of genus 3
is given by
$$
\frac{(p-1)(p-9)(p-11)(p^3 + 20p^2 -349p -3200)}{1451520}.
$$
As for the main term of the number of superspecial curves of genus 3 with
long automorphism of order 2, by Brock \cite{B} it is given by
$$
\frac{(p-1)(p-9)(p^2 -3p -82)}{1152}.
$$
Therefore, we get our asymptotic behavier of the rate.
}

\end{document}